\magnification 1200
\def\R{{\rm I\kern-0.2em R\kern0.2em \kern-0.2em}}
\def\N{{\rm I\kern-0.2em N\kern0.2em \kern-0.2em}}
\def\P{{\rm I\kern-0.2em P\kern0.2em \kern-0.2em}}
\def\B{{\rm I\kern-0.2em B\kern0.2em \kern-0.2em}}
\def\Z{{\rm I\kern-0.2em Z\kern0.2em \kern-0.2em}}
\def\C{{\bf \rm C}\kern-.4em {\vrule height1.4ex width.08em depth-.04ex}\;}

\def\D{{\Delta}}

\def\cH{{\cal H}}

\def\G{{\Gamma}}

\font\ninerm=cmr8
\noindent  {\ninerm To appear in Proc.\ Roy.\ Soc.\ Edinb.\ Sect.\ A} 
\vskip 25mm

\centerline {\bf HOLOMORPHIC EXTENDIBILITY AND MAPPING DEGREE}
\vskip 4mm
\centerline{Josip Globevnik}
\vskip 4mm
{\noindent \ninerm ABSTRACT\ \ Let $D$ be a bounded, finitely connected domain in $\C$ 
without isolated points in the boundary and let $f$ be a continous function on $bD$. 
Let $\tilde f$ be a continuous extension of $f$ to $\overline D$. We prove that $f$ 
extends holomorphically through $D$ if and only if the 
degree of $\tilde f+h$ is nonnegative for every holomorphic function $h$ on $D$ 
such that $\tilde f + h$ is bounded away from zero near $bD$. 
} 
\vskip 6mm

\bf 1.\ Introduction and the main result\rm
\vskip 2mm
H.\ Alexander and J.\ Wermer [AW] obtained a characterization of those compact submanifolds of 
$\C^N$ which are boundaries 
of a\-na\-lytic varieties in terms of their linking numbers with respect to 
algebraic varieties. Their results 
inspired some work on characterizing continuous 
boundary values of holomorphic functions in terms of mapping degree [S] and, in the 
case of one variable, in terms of the argument principle [S, G1, G2] which 
brought some new insights also into the classical one variable theory. 

Let $D$ be a bounded domain in $\C$ and let $A(D)$ be the algebra of
all continuous functions on $\overline D$ which are holomorphic on $D$.

In the special case when $D$ is bounded by finitely many pairwise
disjoint simple closed curves the boundary values 
of functions from $A(D)$ 
can be characterized in terms of the argument principle:
\vskip 2mm
\noindent\bf THEOREM 1.1\rm \ [G2]\ \it Let $D\subset\C$ be a bounded 
domain whose boundary 
consists of a finite 
number of pairwise disjoint simple closed curves. A continuous function 
$f$ on $bD$ extends to a function in $A(D)$ if 
and only if for each $g\in A(D)$ such that $f+g$ has no zero on $bD$, 
the change of argument of $f+g$ along $bD$ is nonnegative.
\vskip 2mm
\rm
We would like to obtain a similar theorem for general domains. Trying to do this we
first formulate Theorem 1.1 in terms of mapping degree.

Let $\Psi$ be a continuous function on $bD$ which does not vanish on $bD$. 
Let $\tilde\Psi$ be a continuous extension of $\Psi$ to $\overline D$. Approximate 
$\tilde \Psi $ on $\overline D$ uniformly by a function $G$ which is 
smooth in an open neighbourhood of $\overline D$ and which does not vanish on $bD$. 
Perturbing $G$ slightly we may assume 
that $0$ is a regular value of $G$ so that $G^{-1}(0)\cap D$ is a finite subset of $D$ 
and each point in 
$G^{-1}(0)\cap D$ is a regular point of $G$. Let $\nu $ be the number of points 
in $G^{-1}(0)\cap D$
at which 
$G$ preserves orientation minus the number of points in  $G^{-1}(0)\cap D$ at which $G$ 
reverses orientation. 
The number $\nu$  depends neither on the choice of the extension $\tilde\Psi $ 
of $\Psi$ nor on the 
choice of $G$ 
provided that $G$ approximates $\tilde\Psi $ on $\overline D$ well enough [D]. We shall 
call the number $\nu $  the \bf degree\rm \ of
$ \Psi $, \ $\nu= deg (\Psi)$ \ .  We mention some of its properties. 
If $\Psi _t,\ 0\leq t\leq 1,$ 
is a continuous family of continuous functions on $bD$ such that $\Psi _t\not = 0$ 
on $bD$ for all $t,\ 0\leq t\leq 1$, then $deg (\Psi _1) = deg (\Psi _0)$. 
If $\tilde\Psi $ is a
continuous function on $\overline D$ and 
$D_1\subset D$ is an open set such that $\tilde \Psi $ has no zero on 
$\overline D\setminus D_1$ 
then $deg(\tilde\Psi |bD) = deg(\tilde\Psi |bD_1).$ Further, if $bD$ 
consists of a finite number 
of pairwise disjoint simple closed curves then $2\pi deg(\Psi)$ equals 
the change of argument 
of $\Psi $ along $bD$\ [D]. In the special case when $\tilde\Psi $ 
is holomorphic on $D$ then $\tilde\Psi 
$ preserves orientation and so the degree of
$\Psi $ equals the number of 
zeros of $\tilde\Psi $ in 
$D$.

We should perhaps point out that usually one calls the number $\nu $ above the degree 
of $\tilde\Psi $ so that one talks about the degree of continuous functions
on $\overline D$ that 
have no zero on $bD$. However, since two such functions have the same degree 
provided that they coincide on $bD$\ [D] one, can, as we do, talk about 
the degree of continuous functions on $bD$ without zeros on $bD$. 

One can rewrite Theorem 1.1. as
\vskip 2mm
\noindent\bf THEOREM 1.2\ \it Let $D$ be as in Theorem 1.1. A continuous 
function $f $ extends holomorphically through 
$D$ if and only if $deg(f + g)\geq 0$  
whenever $g\in A(D)$ is such that 
$f+g$ has no zero on $bD$.
\vskip 2mm
\noindent \bf Example \rm [S]\ Let $\D$ be the open unit disc in $\C$ and 
let $D=\D \setminus 
[0,1)$. Define 
$$
 f(z) =
\left\{
\eqalign{&z\ \ \ \ (z\in b\D)\cr
&1\ \ \ \ (z\in [0,1))\cr}\right. 
$$
The function $f$ is continuous on $bD=b\D\cup [0,1)$. Note that every 
$h\in A(D)$ is actually holomorphic on 
$\D$ so that $A(D)=A(\D)$. Let $h\in A(D)$ be such that $f+h\not= 0 $ on $bD$ 
and let $G$ be a 
continuous extension of $f+h|bD$ to $\overline D = \overline \D.$ Then 
$deg (G|bD) = deg (G|bD_r)$  
where 
$D_r =\D\setminus ([0,1]+r\overline\D )$ provided that $r>0$ is sufficiently small.
The domain $D_r$ is 
bounded by a simple closed curve so the change of argument  of $G$ along $bD_r$ equals 
$2\pi deg(G|bD_r)$. Since $G$ is continuous on $\overline\D$ the change of 
argument of $G$ along 
$\overline D\cap b([0,1]+r\D )$ tends to zero as $r\searrow 0$ and the 
change of argument of $G$ along 
$bD \setminus ([0,1] + r\D )$ tends to the change of argument 
of $G$ along $b\D $ which, 
by the argument principle, equals $2\pi $ times the number of zeros
of $z\mapsto z+h(z)$ in $\D $.
Thus, for any $h\in A(D)$ such that $f+h\not= 0$ on $bD$, the 
degree of $f+h$ is nonnegative, yet $f$ 
does not extend holomorphically through $D$. It follows that in general 
Theorem 1.2 does not hold for more general 
domains than the ones bounded by finitely many pairwise disjoint 
simple closed curves. 

In the present paper we prove an analogue of Theorem 1.2 for  
bounded, finitely connected domains in $\C$ without isolated 
points in the boundary. 

\rm Let $D$ be a bounded domain in $\C$. We first define the degree 
of functions which are continuous and nonzero on $U\cap D$ for some 
open neighbourhood $U$ of $bD$ in $\C$.
\vskip 2mm
\noindent\bf DEFINITION 1.3\ \it Let $U$ be an open neighbourhood of $bD$ in $\C$ 
and let $\Psi $ be a continuous function on $U\cap D$ which has no zero on $U\cap D$. 
Let $W$ be a relatively compact open subset of $D$ which contains $D\setminus U$.  
We define \bf the degree of\ \it
$\Psi $ as 
$deg(\Psi ) = deg(\Psi|bW)).$ 
\vskip 2mm
\noindent \rm By the properties of the degree mentioned in Section 1, \ $deg (\Psi )$ 
is well defined. It 
does not depend on the choice of $W$. Moreover, it does not depend on the 
choice of $U$:\ 
if $U_1\subset U$ is a neighbourhood of $bD$ in $\C$ then $deg(\Psi|U_1\cap D) 
= deg (\Psi )$. 

One can express the degree in terms of the change of the argument. 
Exhaust $D$ by a sequence $D_m$ of domains
$$
D_1\subset\subset D_2\subset\subset\cdots ,\ \cup_{m=1}^\infty D_m=D,
$$
such that for each $m,\ bD_m$ consists of finitely 
many pairwise disjoint simple closed curves. If $\Psi $ is a 
continuous function on $U\cap D$ for some open neighbourhood $U$ of $bD$ in $\C$ such that 
$\Psi\not= 0$ on $U\cap D$ then
there is an $m_0$ such that 
$2\pi deg(\Psi)$ equals the change of argument of $\Psi $ 
along $bD_m$ for all $m\geq m_0$. 
\vskip 2mm\bf \noindent DEFINITION 1.4\ \it Let $\Phi $ be a 
function on $U\cap D$ where U is a neighbourhood of $bD$ in $\C$. 
We say that $\Phi $ is \bf bounded away from zero near\it $\ bD$  if there
are a neighbourhood $V\subset U$ of  $bD$ and $\delta >0$ such that $|\Phi |\geq \delta $ on 
$V\cap D$.\rm
\vskip 2mm
Our main result is
\vskip 2mm
\noindent\bf THEOREM 1.5\ \it Let $D$ be a bounded, finitely connected domain in $\C$ 
without isolated points in the boundary and let $f$ be a continous function on $bD$. 
Let $\tilde f$ be a continuous extension of $f$ to $\overline D$. The function $f$ 
extends holomorphically through $D$ if and only if the degree of $\tilde f + h$ is nonnegative 
for every holomorphic function $h$ on $D$ such that 
$\tilde f+h$ is bounded away from zero near $bD$. 

\vskip 4mm
\bf 2.\ Single valued conjugates
\vskip 2mm \rm 
Let $D$ be as in Theorem 1.5. Exhaust $D$ by a sequence of domains $D_m$ as in 
Section 1. Let $\G_1, \G_2,\cdots \G_n$ be 
the components of $bD$ such that $\Gamma _n$ contains the boundary of the unbounded component of 
$\C\setminus\overline D$.  Observe 
that our conditions imply that $D$ is a Dirichlet
domain, that is, for every continuous function $\phi $ on $bD$ there is a continuous
 extension of $\phi $ to $\overline D$ which is harmonic on $D$. This extension is unique
and we denote it by $\cH (\phi)$. 

There are smoothly bounded Jordan domains $\Omega _j,\ 1\leq j\leq n-1$, such that 
$\Gamma _j\subset \Omega _j\ (1\leq j\leq n-1)$, such that the closures 
$\overline\Omega _j$ are pairwise disjoint and such that 
$b\Omega _j\subset D\ (1\leq j\leq n-1)$. 

A holomorphic function $h$ on $D$ has a single valued primitive function on 
$D$ if and only if
$$
\int_{b\Omega _j}h(z)dz = 0\ \ (1\leq j\leq n-1).
$$
For each $j,\ 1\leq j\leq n-1,$ let $\omega _j$ be the continuous function on $\overline D$ 
which is harmonic on $D$ and which satisfies $\omega _j|\Gamma_j\equiv 1,\ \omega_j|\Gamma _k 
\equiv 0\ (k\not= j,\ 1\leq k\leq n-1)$.  
We shall need the following lemma which is well known for smoothly bounded domains [B].
\vskip 2mm
\noindent\bf LEMMA 2.1\ \it Given a harmonic function $h$ on $D$ there is a 
unique $(n-1)$-tuple $c_1(h),$ $ c_2(h),$ $\cdots c_{n-1}(h)$ of constants such that 
$h+\sum _{j=1}^{m-1}c_j(h)\omega _j $ has a single valued conjugate on $D$. 
\vskip 2mm
\rm For each $k,\ 1\leq k\leq n-1$, let $F_k$ be a multivalued holomorphic function on $D$ 
whose real part is $\omega _k$. Then $F_k^\prime $ is single valued on $D$ and 
for each $j,\ 1\leq j\leq n-1,$ the integral
$$
\int _{b\Omega _j} F_k^\prime (z)dz =2\pi i \alpha_{k,j}
$$
is purely imaginary. We use the reasoning from [E, p.276] to 
show that the (real) matrix $[\alpha_{k,j}]$ is nonsingular. Suppose that this 
is not the case. 
Then there are real constants $\lambda _1, \lambda _2,\cdots \lambda _{n-1}$,
not all zero, such that
$$
\sum _{k=1}^{n-1}\lambda _k\alpha _{k,j} = 0\ \ (1\leq j\leq n-1)
$$
which implies that
$$
\int_{b\Omega_j}\bigl[\sum_{k=1}^{n-1}\lambda _kF_k(z)\bigr] dz = 0
\ \ \ (1\leq j\leq n-1)
$$
so the function $F =\sum _{k=1}^{n-1}\lambda _kF_k$ is 
a single valued holomorphic function on $D$. Note that the real part $\Re F$ 
has a continuous extension 
$\tilde{\Re F} $ to $\overline D$ and
$$
\tilde{\Re F}|\Gamma_k\equiv \lambda _k\ (1\leq k\leq n-1),
\ \tilde{\Re F}|\Gamma_n\equiv 0.
$$
Suppose that $F$ takes in $D$ a value $w_0$ that does not belong
to the union of vertical lines
$i\R,\ \lambda_1+i\R,\ \lambda_2+i\R,\cdots , \lambda _{n-1}+i\R$.
Choose $\delta >0$ so small 
that $w_0$ does not belong to the union of vertical strips 
$(-\delta, \delta)+i\R,$ $ 
(\lambda_1-\delta, \lambda_1+\delta)+i\R,\ \cdots\ (\lambda_{n-1}
- \delta,
\lambda_{n-1}+\delta )+i\R$. Let $w_0=F(z_0)$ where
$z_0\in D$ and notice that $m$ can be chosen so 
large that $z_0\in D_m$ and that $F(bD_m)$ is contained 
in the union of the vertical strips above. This 
implies that the change of argument of $z\mapsto F(z)-w_0$ 
along $bD_m$ is zero, contradicting the fact 
that $F(z_0)=w_0$ for $z_0\in D_m$. Thus $F(D)$ is 
contained in the union of the vertical lines above and 
consequently $F$ is a constant which is not possible since
at least one of $\lambda_k's$ is different from 
0, say, $\lambda_j\not=0$ and 
$\tilde{\Re F}|\Gamma_j\equiv\lambda_j,\ \tilde{\Re F}|\Gamma_n\equiv 0$. 
The contradiction shows that the matrix  $[\alpha_{k,j}]$ is nonsingular.
\vskip 2mm\noindent
\bf Proof of Lemma 2.1.\ \rm It suffices to prove that given
a real harmonic function $h$ on $D$ there are 
unique real constants $c_1(h),\cdots ,c_{n-1}(h)$ such that
$h+\sum_{j=1}^{n-1}c_j(h)\omega_j$ has 
a single valued conjugate. This is the same as 
to say that there are unique real constants $c_1(h),\cdots ,c_{n-1}(h)$ such that 
$ H+\sum_{j=1}^{n-1}c_j(h)F_j$ is single valued where $H$ is a multiple
valued holomorphic function on $D$ whose real part is $h$. This happens if and only if
$$
\int _{b\Omega_j}\bigl[H^\prime (z)+\sum_{k=1}^{n-1}c_k(h)F_k^\prime (z)\bigr]dz=0\ \ 
(1\leq j\leq n-1)
$$
so
$$
\sum_{k=1}^{n-1}c_k(h)\bigl[{1\over{2\pi i}}\int_{b\Omega_j}F_k^\prime (z)dz\bigr] =
-{1\over{2\pi i}}\int_{b\Omega _j}H^\prime (z)dz \ \ (1\leq j\leq n-1).
\eqno (2.1)
$$
Since the matrix $[\alpha_{j,k}]$ is nonsingular there is a unique  $(n-1)$-tuple 
$c_1(h),\cdots,c_{n-1}(h)$ of real constants such that (2.1) is satisfied. 
This completes the proof.
\vskip 4mm
\bf 3.\ Proof of Theorem 1.5\rm
\vskip 2mm
Observe first that if $U$ is a neighbourhood 
of $bD$ in $\C$ and  $\Phi,\ \Psi$ are continuous functions on
$U\cap \overline D$ such that $\Phi\equiv \Psi $ on $bD$, and $h$ 
is a continuous 
function on $U\cap D$ such that $\Phi+h$ is bounded away from $0$ 
near $bD$,  
then $\Psi + h$ is bounded away from $0$ near $bD$ and \ 
$deg(\Phi+h) = deg(\Psi +h)$. To see this notice first that by the assumption there are a 
neighbourhood $V\subset U$ of $bD$ and 
$\delta >0$ such that $|\Phi+h|\geq \delta $ on $V\cap D$. 
 Since $\Phi \equiv\Psi $ on
$bD$ and since $\Phi$ and $\Psi $ are continuous there is a neighbourhood $V_1\subset V$ 
of $bD$ in $\C$ such that $|\Phi -
\Psi|<\delta /2$ on $V_1\cap D$. For each $t,\ 0\leq t\leq 1$, we have 
$|\Phi +t(\Psi-\Phi)+h|\geq|\Phi+h|-\delta/2\geq \delta/2$ 
on $V_1$, that is, 
for each $t,\ 0\leq t\leq 1$, the function $(1-t)(\Phi+h)+t(\Psi +h)$ 
is bounded 
away from $0$ on $D\cap V_1$. In particular, $\Psi + h$ is bounded 
away from $0$ on $D\cap V_1$. Moreover, since $(\Psi+h)|V_1$ is homotopic to $(\Phi +h)|V_1$ 
through a family of maps with ranges contained in $\C\setminus\{ 0\}$ it follows that 
$deg(\Psi +h)=deg (\Phi +h)$. 

We now turn to the proof of Theorem 1.5.  Denote by
$Z$ the identity, $Z(z)\equiv z$.
Suppose that $f$ is a continuous 
function on $bD$ which has a continuous extension $g$ to $\overline D$ which 
is holomorphic on $D$. Let $h$ be a holomorphic function on $D$ such that $g+h$ is 
bounded away from zero near $bD$. This means that there are an open 
neighbourhood $U$ of $bD$ in $\C$ and a $\delta >0$ such that $|g+h|\geq \delta>0$ on $U\cap D$.
Choose $m$ so large that $bD_m\subset U$. Then $2\pi deg (g+h)$ equals 
the change of 
argument of $g + h$ along $bD_m$. Since $g$ and $h$ are holomorphic on $D$ 
the argument principle implies that the change  
of argument of $g+h$ along $bD_m$ equals $2\pi$ times the number 
of zeros of $g +h$ 
in $D_m$ so it is nonnegative. Thus $deg (g+h)\geq 0$. The resoning above  implies 
that $\tilde f+h$ is
bounded away from zero near $bD$ and that $deg (\tilde f+h) = deg (g+h)\geq 0$.
This proves the only if 
part of the theorem.

The proof of the if part is very short in the special case when $D$ is simply connected.
It uses the basic fact that a continuous function $f$ on 
$bD$ extends holomorphically through  $D$ if and only if the harmonic extension of $Zf$ 
equals $z$ times the harmonic extension of $f$, that is, if and only 
if $\cH (Zf)(z) \equiv z\cH(f)(z)\ \ (z\in D)$ [G2]. So suppose that $f$ does  not extend 
holomorphically through $D$. Then 
 there is an $a\in D$ such that $\cH ((Z-a)f)(a)\not=0$, and multiplying $f$ by $
 e^{i\gamma},\ \gamma\in\R$, we may assume that $\cH ((Z-a)f)(a)=\alpha $ where $\alpha>0$. 
 Since $D$ is simply connected it follows that $\cH ((Z-a)f)(z) - \alpha = (z-a)F(z) - 
 \overline{(z-a)G(z)}$ where $F$ and $G$ are holomorphic on $D$ and so 
 $\cH ((Z-a)f)(z) - (z-a)F(z) - 
(z-a)G(z) \in\alpha+i\R \ (z\in D)$ which implies that $z\mapsto (1/(z-a))\cH ((Z-a)f)(z)-
F(z)-G(z)$ is bounded away from zero near $bD$ and its degree equals $-1$.

The proof of the if part in the case when $D$ is multiply connected is more complicated 
because in general, conjugate functions are no more single valued. 
We shall use the reasoning from [G2]. Assume that $f$ is a continuous 
function on $bD$ that does
not extend holomorphically through $D$. Denote by
$Z$ the identity, $Z(z)\equiv z$ and define 
$$
A(a,f)=\cH ((Z-a)f)(a) = \cH (Zf)(a)- a\cH (f)(a)\ \ (a\in\overline D).
$$
Since $\cH (f)$ is not holomorphic on $D$ it follows that $\{ a\in D\colon \ A(a,f)=0\} $ 
is a closed, nowhere dense subset of $D$ [G2]. There 
are constants $c_k,\ d_k, \ 1\leq k\leq n-1$, 
such that
$$
\cH (f)(z) +\sum_{k=1}^{n-1}c_k \omega _k(z) =
P(z)+\overline{Q(z)}\ \ (z\in D)
$$
$$
\cH (Zf)(z) +\sum_{k=1}^{n-1}d_k \omega _k(z) =
R(z)+\overline{S(z)}\ \ (z\in D)
$$
where $P, Q, R, S$ are holomorphic functions on $D$. Define
$$
\Phi_a(z)= \sum_{j=1}^{n-1} [d_j-ac_j].[\omega_j(z)-\omega _j(a)] - A(a,f).
$$
For each $a\in D$ the function $\Phi _a$ is continuous on $\overline D$ 
and harmonic 
on $D$ and by the preceding 
discussion, the function $z\mapsto \cH((Z-a)f)(z)+\Phi _a(z)\ \ (z\in D)$ 
vanishes at $a$ and 
has a conjugate on $D$, that is, it is of the form 
$(z-a)F_a(z)+\overline{(z-a)G_a(z)}$ where 
$F_a$ and $G_a$ are single valued holomorphic functions on $D$. 

The function $\Phi _a$ is constant on each $\Gamma _k, 1\leq k\leq n$. 
Repeating the proof of Lemma 6.1 in [G2] we see that there 
is an $a\in D$ such that all these constants are different from $0$, 
that is, $\Phi _a(z)\not= 0 \ (z\in bD)$. Replacing $f$ by $e^{i\gamma}f,\ \gamma\in\R$ 
will multiply all these constants with $e^{i\gamma}$ so we may assume with no 
loss of generality that the real parts of all these constants are different from zero, 
that is,  $\Re \Phi_a (z)\not= 0\ (z\in bD)$. Hence, if 
$$
W(z)=\cH ((Z-a)f)(z) - (z-a)F_a(z)-(z-a)G_a(z)\ \ (z\in D)
$$
it follows that 
$$
\Re W(z) = - \Re \Phi _a(z)\ \ (z\in D)
$$
where the function on the right is continuous on $\overline D$ and for 
each $j,\ 1\leq j\leq n,$ equal to 
a nonzero real constant $\beta _j$ on $\Gamma _j$. 
This shows that $W$ is bounded away from zero on $U\cap D$ for a 
neighbourhood $U$ of $bD$.         
Choose $\delta >0$ so small that $0\not\in \cup_{j=1}^n [\beta_j-\delta,
\beta _j+\delta]$.      
There is an $m_0\in \N$ such that for all $m\geq m_0,$ we have 
$W(bD_m)\subset \cup_{j=1}^n 
[\beta_j-\delta, \beta_k+\delta]+i\R$. It follows that the change 
of argument of $W$ along $bD_m$ is zero,
and consequently, provided that $a\in D_{m_0}$, the change of
argument of $z\mapsto W(z)/(z-a)$
along $bD_m$ equals $-2\pi$ for all $m\geq m_0$. This implies that the function 
$$
z\mapsto {1\over{z-a}}W(z) = {1\over{z-a}}\cH ((Z-a)f)(z) -F_a(z)-G_a(z)
$$
is bounded away from zero near $bD$ and its degree equals $-1$. Since the function 
$z\mapsto (1/(z-a))\cH ((Z-a)f)(z)$ is a continuous extension of $f$ to 
$\overline D\setminus\{ 0\}$ 
Proposition 3.1 implies that $\tilde f-F_a-G_a$ is bounded away from zero near $bD$ and \ 
$deg (\tilde f-F_a-G_a)=-1$.
This completes the proof.
\vskip 4mm
\bf 4.\ An example \rm
\vskip 2mm
Let $\D $ be the open unit disc in $\C$ and let $D=\D\setminus\{ 0\}$. Define
$$
f(z) = \left\{\eqalign{\ &0 \ \ (z\in b\D)\cr
\ &1\ \ (z=0).\cr}\right.
$$
Suppose that $\tilde f$ is a continuous extension of $f$ to $\overline D =\overline\D $ and 
let $h$ be a holomorphic function on $D$  such that $\tilde f+h$ 
is bounded away from 0 on $bD$. In particular, there are $\delta >0$, a 
neighbourhood $U$ of $b\D$ in $\C $ and a neighbourhood $V$ of $0$ such that
$|\tilde f + h|\geq 2\delta$ on $U\cap \D$ and on $V\setminus \{ 0\}$. 
Shrinking $U$ and $V$ if necessary we may assume that 
$$
|h|\geq \delta\  \hbox{\ on\ } U\cap\D
\eqno (4.1)
$$
$$
|h+1|\geq \delta\ \hbox{\ on\ } V\setminus\{ 0\}.
\eqno (4.2)
$$
We have  
$2\pi deg(\tilde f+h) = V_R-V_\rho$ where $V_R$ is the change  
of argument of $\tilde f+h$ along $b(R\D )$ for $R<1$ very 
close to $1$ and $V_\rho $ is the change of argument 
of $\tilde f + h $ along $b(\rho\D )$ for $\rho >0$ very small.
Shrinking $U$ if necessary we may assume that
$|\tilde f|<\delta/2$ on $U\cap \D $. Thus, if $R<1$ is sufficiently close to $1$,
\ $t\mapsto (h+t\tilde f)|b(R\D )\ (0\leq t\leq 1)$ is a homotopy between $h|b(R\D)$ and
$(h+\tilde f)|b(R\D)$ in $\C\setminus\{ 0\}$ so $V_R$ equals 
the change of argument of $h$ along $b(R\D )$. Similarly, shrinking 
$V$ if necessary we may assume that $|\tilde f
-1|<\delta/2$ on $V\setminus\{ 0\}$, so, provided that
$\rho>0$ is sufficiently small, 
$t\mapsto [h+1 + t(\tilde f-1)]|b(\rho\D)\ (0\leq t\leq 1)$ is a homotopy between 
$(h+1)|b(\rho\D)$ and $ [h+\tilde f]|b(\rho\D)$ in $\C\setminus \{ 0\}$ so $V_\rho$ equals the 
change of argument of $h+1$ along $b(\rho\D)$. The condition (4.2) implies that $h$ is either 
holomorphic at $0$ or has a pole at $0$. In the first case the change of argument of $h+1$ 
along $b(\rho\D)$ 
is zero since $h+1$ has no zero on $\rho\overline\D$. In the second case there is a $\rho >0$
such that $h$ is large on $b(\rho\D)$, say, $|h|>10$. This implies that the change of argument of $h+1$ along 
$b(\rho\D)$ equals the change of argument of $h$ along $b(\rho\D)$. Thus, in the first case, $V_\rho =0$ and so $V_R-V_\rho $ equals 
the number of zeros of $h$ on $R\D $ and so is nonnegative. In the second case $V_R-V_\rho $ 
equals the number of zeros of $h$ on $RD \setminus \rho\overline D$ and so is nonnegative. Thus, $deg(f+h)\geq 0$ 
for every holomorphic function $h$ on $D$ such that $f+h$ is bounded away from $0$ on $bD$, yet $f$
does not extend holomorphically through $D$. This shows that in Theorem 1.3 the assumption that all components of $bD$ 
be nondegenerate cannot be omitted.
\vskip 5mm
This work was supported 
in part by the Ministry of Higher Education, Science and Technology of Slovenia  
through the research program Analysis and Geometry, Contract No.\ P1-0291. 
\vfill
\eject
\centerline{\bf REFERENCES}
\vskip 5mm
\vskip 5mm
\noindent [AW]\ H.\ Alexander and J.\ Wermer: Linking 
numbers and boundaries of varieties. 

\noindent Ann.\ Math.\ 151 (2000) 125-150
\vskip 2mm
\noindent [B]\ S.\ Bell: \it The Cauchy transform, Potential Theory, and 
Conformal Mapping. \rm CRC Press, 
Boca Raton, 1992
\vskip 2mm
\noindent [D]\ K.\ Deimling:\ \it Nonlinear functional analysis.\rm \ Springer Verlag, Berlin,  1980
\vskip 2mm
\noindent [E]\ M.\ A.\ Evgrafov:\ \it Analytic functions.\rm\ Saunders, 
Philadelphia and London, 1966
\vskip 2mm
\noindent [G1]\ J.\ Globevnik:\ Holomorphic extendibility and the argument principle.

\noindent Complex Analysis and Dynamical Systems II. Contemp.\ Math.\ 382 (2005) 171-175
\vskip 2mm
\noindent [G2]\ J.\ Globevnik:\ The argument principle and holomorphic extendibility.

\noindent Journ.\ d'Analyse.\ Math.\ 94 (2004) 385-395
\vskip 2mm
\noindent [S]\ E.L.Stout: Boundary values and mapping degree.

\noindent Michig.\ Math.\ J.\ 47 (2000) 353-368
\vskip 2mm
\noindent [T]\ M.\ Tsuji:\ \it Potential theory in modern function theory.  \rm
Maruzen, Tokyo, 1959
\vskip 10mm
\noindent Institute of Mathematics, Physics and Mechanics

\noindent University of Ljubljana, Ljubljana, Slovenia

\noindent josip.globevnik@fmf.uni-lj.si

\bye